\newfont{\Bbb}{msbm10 scaled\magstephalf}
\documentclass[reqno]{amsart}
\usepackage{amssymb}
\usepackage{amsfonts}
\usepackage{amsmath, amsthm, amscd, amsfonts}
\usepackage{cases}

\numberwithin{equation}{section}

\begin{document}
\title{Difference of Weighted Composition Operators from $\alpha$-Bloch spaces to $\beta$-Bloch spaces}

\author[N. Xu and Z.H. Zhou]  {Ning Xu and Ze-Hua Zhou$^*$}

\address{\newline Ning Xu\newline School of science, Jiangsu ocean University\newline Lianyungang 222005\newline P.R. China}
\email{gx899200@126.com}

\address{\newline  Ze-Hua Zhou\newline  School of Mathematics, Tianjin University, Tianjin 300354, P.R. China.\newline
} \email{zehuazhoumath@aliyun.com;
zhzhou@tju.edu.cn}



\date{}
\thanks{\noindent $^*$Corresponding author.\\
The work was supported in part by the National Natural Science Foundation of
China(Grant No.11771323) and the Scientific Research Foundation for Ph.D.of Jiangsu Ocean University(No.KQ17006).}

\begin{abstract}
In this paper, we study the boundedness and compactness of the differences of two weighted composition operators acting from $\alpha$-Bloch space to $\beta$-Bloch space on the open unit disk. This study has a relationship to the topological structure of weighted composition from $\alpha$-Bloch space to $\beta$-Bloch space.

\end{abstract}

\maketitle

\section{Introduction}
Let $\mathbb{D}$ be an open disk in the complex plane $\mathbb{C}$ and $H(\mathbb{D})$ be the class of all functions analytic in $\mathbb{D}$.
We denote by $S(\mathbb{D})$ the set of all analytic self-maps on $\mathbb{D}$. For $u\in H(\mathbb{D})$, denote by $M_u$ the multiplication operator.
For $\varphi\in S(\mathbb{D})$, denote by $C_\varphi$ the composition operator. Then the weighted composition operator $uC_\varphi$ is a linear transformation of $H(\mathbb{D})$ denote by
\begin {eqnarray*}
(uC_\varphi)f(z)=(M_uC_\varphi)f(z)=u(z)f(\varphi(z)).
\end {eqnarray*}
For $0<\alpha<\infty$, recall that the Bloch type space $\mathcal{B}^\alpha$, or $\alpha$-Bloch space, consists of all $f\in H(\mathbb{D})$ such that
\begin{eqnarray*}
\|f\|_{\alpha}=\sup\limits_{z\in\mathbb{D}}(1-|z|^2)^\alpha|f'(z)|<\infty.
\end{eqnarray*}
It is well known that $\mathcal{B}^\alpha$ is a Banach space under the norm $\|f\|_{\mathcal{B}^\alpha}=|f(0)|+\|f\|_\alpha$. The little $\alpha$-Bloch $\mathcal{B}_0^\alpha$
is a subspace of $\mathcal{B}^\alpha$, consisting of all $f\in H(\mathbb{D})$ such that $\lim\limits_{|z|\rightarrow 1}(1-|z|^2)^\alpha|f'(z)|=0$.
When $\alpha=1$, $\mathcal{B}^\alpha$ is classical space $\mathcal{B}$.
\

For $z,w\in \mathbb{D}$, the pseudo-hyperbolic distance $\rho(z,w)$ between $z$ and $w$ is given by
\begin {eqnarray*}
\rho(z,w)=\Big|\frac{w-z}{1-\overline{w}z}\Big|.
\end{eqnarray*}
For $a\in \mathbb{D}$, Let $\sigma_a$ be the M\"{o}bius transformation of $\mathbb{D}$ defines by
\begin{eqnarray*}
\sigma_a(z)=\frac{a-z}{1-\overline{a}z}.
\end{eqnarray*}
We remark that $\rho(a,z)=|\sigma_a(z)|\leq 1$. For $\varphi\in S(\mathbb{D})$, the Schwarz-Pick type derivative of $\varphi$ is defined by
\begin{eqnarray*}
\varphi^\sharp(z)=\frac{(1-|z|^2)^\beta}{(1-|\varphi(z)|^2)^\alpha}\varphi'(z).
\end{eqnarray*}
The Schwarz-Pick Lemma implies that $\|\varphi^\sharp\|_\infty\leq 1$. Madigan and Matheson \cite{MM} proved that $C_\varphi$ is compact on $\mathcal{B}$ if and only if $|\varphi^\sharp(z)|\rightarrow 0$ wherever $|\varphi(z)|\rightarrow 1$. Ohno and Zhao \cite{OZ} generalized this result to the case of $uC_\varphi$
on  $\mathcal{B}$.
\

The study of the difference of two composition operators was started on the Hardy space $H^2$. The main purpose for this study is to understand the topological structure of $\mathcal{C}(H^2)$, the set of composition operators on $\mathcal{C}(H^2)$, see \cite{SS}. Then, MacCluer, Ohno and Zhao \cite{MOZ} considered the above problems on H$^\infty$. This work gave a relationship between a component problem and the behavior of the difference of two composition operators acting from $\mathcal{B}$ to H$^\infty$. After that, such related problems have been studied on several spaces of analytic functions by many authors, see \cite{HO1,FZ1,KSK,SL} for the study of differences of composition operators on various spaces. Some related results concerning the differences of weighted composition operators on various spaces can be founded in, for example, \cite{CZ,FZ2,H,HO2}.
\

As an analogue, the topological structure of the set $\mathcal{C}(\mathcal{B}^\alpha\rightarrow\mathcal{B}^\beta)$ of bounded weighted composition operators from $\mathcal{B}^\alpha$ to $\mathcal{B}^\beta$ can be considered. In this context, we deal with the differences of weighted composition operators from $\mathcal{B}^\alpha$ to $\mathcal{B}^\beta$ when $\alpha>1$. The main purpose of this paper is to express the boundedness and compactness of  $uC_\varphi-vC_\psi$ from $\mathcal{B}^\alpha$ to $\mathcal{B}^\beta$ which generalizes \cite{H}. The authors expect that these results will play some roles in the succeeding investigation.
\

For two quantities $A$ and $B$ which may depend on $\varphi$ and $\psi$, we use the abbreviation $A\lesssim B$ whenever this is a positive constant $C$(independent of $\varphi$ and $\psi$) such that $A\leq CB$. We write $A\sim B$ if $A\lesssim B \lesssim A$.

\section{Prerequisites}
We collect some basic properties of functions in $\mathcal{B}^\alpha$. It is known that the follows hold(see \cite{Z}), for $z,w\in \mathbb{D}$
\begin{eqnarray*}
d_\alpha(z,w)=\sup\{|f(z)-f(w)|:\|f\|_{\mathcal{B}^\alpha}\leq 1\},
\end{eqnarray*}
\begin{eqnarray*}
|f(z)|\leq\frac{\|f\|_{\mathcal{B}^\alpha}}{(1-|z|^2)^{\alpha-1}},\ \  f\in\mathcal{B}^\alpha, \ \ \alpha>1.
\end{eqnarray*}
Moreover we use the following notion:
\begin{eqnarray*}
\rho(z)=\rho(\varphi(z),\psi(z)),\ \ \ \sigma(z)=\frac{\varphi(z)-\psi(z)}{1-\overline{\varphi(z)}\psi(z)},
\end{eqnarray*}
\begin{eqnarray*}
\tau(z)=\frac{(1-|\varphi(z)|^2)(1-|\psi(z)|^2)}{(1-\overline{\varphi(z)}\psi(z))^2}.
\end{eqnarray*}
We remark that
\begin{eqnarray*}
|\sigma_{\varphi(w)}(\psi(w))|=|\sigma(w)|=\rho(w),\ \ \rm{and} \ \ |\tau(z)|=1-\rho^2(z).
\end{eqnarray*}
For $\{z_n\}\subset \mathbb{D}$, denote that $\varphi_n=\varphi(z_n)$ and $\psi_n=\psi(z_n)$. Similarly we use the notion $\sigma_n, \tau_n$.

Note that
\begin{eqnarray*}
U(z)=(1-|z|^2)^\beta u'(z),\ \ \ V(z)=(1-|z|^2)^\beta v'(z).
\end{eqnarray*}

{\bf Lemma 2.1.{\cite{SL}}} Let $0<\alpha<\infty$. For all $z,w\in \mathbb{D}$, the Bloch-type induced distance is given by
\begin{eqnarray*}
\flat_\alpha(z,w)=\sup\limits_{\|f\|_{\mathcal{B}^\alpha}\leq 1}|(1-|z|^2)^\alpha f'(z)-(1-|w|^2)^\alpha f'(w)|.
\end{eqnarray*}
Then
\begin{eqnarray*}
\flat_\alpha(z,w)\lesssim \rho(z,w).
\end{eqnarray*}

{\bf Remark 2.1.} Especially, if we take the function $f_\lambda(z)=\frac{(1-|\lambda|^2)^{\alpha-1}}{(1-\overline\lambda z)^{2(\alpha-1)}}$, $h_\lambda(z)=\int^z_0f_\lambda(\zeta)d\zeta$, then $\|h_\lambda\|_{\mathcal{B}^\alpha}\leq 1$. From Lemma 2.1, take $\lambda=\varphi(w)$, we have that
\begin{eqnarray*}
|1-\tau^{\alpha-1}(w)|&=&\Big|1-\frac{(1-|\varphi(w)|^2)^{\alpha-1}(1-|\psi(w)|^2)^{\alpha-1}}{(1-\overline{\varphi(w)}\psi(w))^{2(\alpha-1)}}\Big|\\
&=&|(1-|\varphi(w)|^2)^{\alpha-1}h'_{\varphi(w)}(\varphi(w))-(1-|\psi(w)|^2)^{\alpha-1}h'_{\varphi(w)}(\psi(w))|\\
&\leq&\flat_\alpha(\varphi(w),\psi(w))\lesssim \rho(w).
\end{eqnarray*}
The result is right for $|1-\tau^{\alpha}(w)|$ too.

{\bf Lemma 2.2.} There exists a constant $C>0$ such that for any $a,b\in  \mathbb{D}$,
\begin{eqnarray*}
\Big|\frac{1-|a|^2}{1-\overline{a}b}\Big|^{\alpha-1}\leq C.
\end{eqnarray*}

{\bf Proof.} It is very easy and we omit.
\section{The boundedness from $\mathcal{B}^\alpha$ to $\mathcal{B}^\beta$}
Here we formulate and prove the main results of this paper.

{\bf Proposition 3.1.}
Let $\varphi,\psi\in S(\mathbb{D})$ and $u,v\in H(\mathbb{D})$.
Then $uC_\varphi-vC_\psi$ is bounded from $\mathcal{B}^\alpha(\alpha>1)$ to $\mathcal{B}^\beta$, if
the following conditions hold:

(i)
$
\sup\limits_{z\in\mathbb{D}}|u(z)\varphi^\sharp(z)-v(z)\psi^\sharp(z)|<\infty.
$

(ii)
(a)
$
\sup\limits_{z\in\mathbb{D}}|u(z)\varphi^\sharp(z)|\rho(z)<\infty,
$
(b)
$
\sup\limits_{z\in\mathbb{D}}|v(z)\psi^\sharp(z)|\rho(z)<\infty.
$

(iii)
$
\sup\limits_{z\in\mathbb{D}}L(z)<\infty,
$\\
where $L(z)=\min\Big\{\frac{|U(z)-V(z)|}{(1-|\varphi(z)|^2)^{\alpha-1}},\frac{|U(z)-V(z)|}{(1-|\psi(z)|^2)^{\alpha-1}},
|U(z)|d_\alpha(z),|V(z)|d_\alpha(z)\Big\}$.

{\bf Proof.} For any $f\in \mathcal{B}^\alpha$ with $\|f\|_{\mathcal{B}^\alpha}\leq 1$, we have that
\begin {eqnarray*}
\|(uC_\varphi-vC_\psi)(f)\|_{\mathcal{B}^\beta}
&=&|u(0)f(\varphi(0))-v(0)f(\psi(0))|+\|(uC_\varphi-vC_\psi)(f)\|_\beta\\
&\leq&|u(0)-v(0)|f(\varphi(0))+|v(0)||f(\varphi(0))-f(\psi(0))|\\
&+&\|(uC_\varphi-vC_\psi)(f)\|_\beta
\end {eqnarray*}
and
\begin {eqnarray*}
\|(uC_\varphi-vC_\psi)(f)\|_\beta
&\leq&\sup_{z\in\mathbb{D}}\Big|u(z)\frac{(1-|z|^2)^\beta}{(1-|\varphi(z)|^2)^\alpha}\varphi'(z)(1-|\varphi(z)|^2)^\alpha f'(\varphi(z))\\
& &-v(z)\frac{(1-|z|^2)^\beta}{(1-|\psi(z)|^2)^\alpha}\psi'(z)(1-|\psi(z)|^2)^\alpha f'(\psi(z))\Big|\\
& &+\sup_{z\in\mathbb{D}}|(1-|z|^2)^\beta u'(z)f(\varphi(z))-(1-|z|^2)^\beta v'(z)f(\psi(z))|\\
&\triangleq&I+II.
\end {eqnarray*}
First we estimate $I$,
\begin {eqnarray*}
I&\leq&\sup\limits_{z\in\mathbb{D}}|u(z)\varphi^\sharp(z)-v(z)\psi^\sharp(z)|\|f\|_{\mathcal{B}^\alpha}
+\sup_{z\in\mathbb{D}}|v(z)\psi^\sharp(z)|\flat_\alpha(\varphi(z),\psi(z))\\
&\leq&\sup\limits_{z\in\mathbb{D}}|u(z)\varphi^\sharp(z)-v(z)\psi^\sharp(z)|+C\sup_{z\in\mathbb{D}}|v(z)\psi^\sharp(z)|\rho(z).
\end {eqnarray*}
Similarly
$I\leq\sup\limits_{z\in\mathbb{D}}|u(z)\varphi^\sharp(z)-v(z)\psi^\sharp(z)|+C\sup\limits_{z\in\mathbb{D}}|u(z)\varphi^\sharp(z)|\rho(z).
$\\
Then we estimate $II$
\begin {eqnarray*}
II&\leq&\sup\limits_{z\in\mathbb{D}}\Big[|(U(z)-V(z))f(\varphi(z))|+|V(z)(f(\varphi(z))-f(\psi(z)))|\Big]\\
&\leq&\sup\limits_{z\in\mathbb{D}}\Big[\frac{|U(z)-V(z)|}{(1-|\varphi(z)|^2)^{\alpha-1}}+|V(z)d_\alpha(z)|\Big].
\end {eqnarray*}
Similarly $II\leq\sup\limits_{z\in\mathbb{D}}\Big[\frac{|U(z)-V(z)|}{(1-|\psi(z)|^2)^{\alpha-1}}+|U(z)d_\alpha(z)|\Big].
$\\
Note
$$
L(z)=\min\Big\{\frac{|U(z)-V(z)|}{(1-|\varphi(z)|^2)^{\alpha-1}},\frac{|U(z)-V(z)|}{(1-|\psi(z)|^2)^{\alpha-1}},
|U(z)d_\alpha(z)|,|V(z)d_\alpha(z)|\Big\}.
$$
Consider $I$ and $II$, then the proposition holds.

{\bf Proposition 3.2.}
Let $\varphi,\psi\in S(\mathbb{D})$ and $u,v\in H(\mathbb{D})$.
If $uC_\varphi-vC_\psi$ is bounded from from $\mathcal{B}^\alpha(\alpha>1)$ to $\mathcal{B}^\beta$, then
the following conditions hold:

(i)
(a)
$
\sup\limits_{z\in\mathbb{D}}|u(z)\varphi^\sharp(z)|\rho^2(z)<\infty,
$
(b)
$
\sup\limits_{z\in\mathbb{D}}|v(z)\psi^\sharp(z)|\rho^2(z)<\infty.
$

(ii)
(a)
$
\sup\limits_{z\in\mathbb{D}}|U(z)|\rho^3(z)<\infty,
$
(b)
$
\sup\limits_{z\in\mathbb{D}}|V(z)|\rho^3(z)<\infty.
$

(iii)
$
\sup\limits_{z\in\mathbb{D}}|u(z)\varphi^\sharp(z)+v(z)\psi^\sharp(z)|\rho(z)<\infty.
$

{\bf Proof.}
(i) Suppose that $uC_\varphi-vC_\psi$ is bounded from $\mathcal{B}^\alpha$ to $\mathcal{B}^\beta$.
Fix $w\in\mathbb{D}$ and take $\lambda=\varphi(w),\mu=\psi(w), z=w$. We define the function 
$
f(z)=f_\lambda(z)\sigma_\lambda(z)\sigma_\mu^2(z).
$
It is easy to see that $f\in \mathcal{B}^\alpha$. By estimate on $\|(uC_\varphi-vC_\psi)(f)\|_\beta$, then we have that
\begin {eqnarray*}
C&\geq&(1-|w|^2)^\beta|u(w)\frac{1}{(1-|\varphi(w)|^2)^{\alpha}}\varphi'(w)\rho^2(\varphi(w),\psi(w))|\\
&=&|u(w)\varphi^\sharp(w)|\rho^2(w).
\end {eqnarray*}
Similarly define
$
f(z)=f_\mu(z)\sigma_\mu(z)\sigma_\lambda^2(z),
$
we obtain
$
|v(w)\psi^\sharp(w)|\rho^2(w)\leq C.
$
Since $w$ is arbitrary in $\mathbb{D}$, we get (i).

(ii) Consider the function $g(z)=f_\mu(z)\sigma^3_\lambda(z),$ it is easy to see that $g(z)\in \mathcal{B}^\alpha$.
Take $\lambda=\varphi(w),\mu=\psi(w), z=w$,
and by estimate on $\|(uC_\varphi-vC_\psi)(g)\|_\beta$, we have that
\begin {eqnarray*}
C&\geq&|(1-|w|^2)^\beta\Big|v'(w)\frac{\rho^3(w)}{(1-|\psi(w)|^2)^{\alpha-1}}\Big|\\
& &-|v(w)\psi^\sharp(w)|\rho^2(w)\Big|2(\alpha-1)\overline{\psi(w)}\sigma(w)+3(1-|\psi(w)|^2)\sigma'_{\varphi(w)}(\psi(w))\Big|\\
&\geq&|V(w)|\frac{\rho^3(w)}{(1-|\psi(w)|^2)^{\alpha-1}}
-C|v(w)\psi^\sharp(w)|\rho^2(w)
\end {eqnarray*}
where we use the fact that
\begin {eqnarray*}
|2(\alpha-1)|\overline{\psi(w)}\sigma_{\varphi(w)}(\psi(w))-3\tau(w)|\leq C,\ \
\rm{and} \ \ |v(w)\psi^\sharp(w)|\rho^2(w)\leq C,
\end {eqnarray*}
therefore we have that
$$
|V(w)|\rho^3(w)\leq C(1-|\psi(w)|^2)^{\alpha-1}\leq C.
$$
Similarly define
$
g(z)=f_\lambda(z)\sigma^3_\mu(z),
$
we get
$
|U(w)|\rho^3(w)\leq C.
$
Since $w$ is arbitrary in $\mathbb{D}$, we get (ii).

(iii) Fix $w\in \mathbb{D}$ and define $h(z)=f_\mu(z)\sigma_\lambda(z)\sigma_\mu(z)$, then $ h\in \mathcal{B}^\alpha$.
Take $\lambda=\varphi(w),\mu=\psi(w), z=w$, we obtain that
\begin {eqnarray*}
C&\geq&\Big|u(w)\varphi^\sharp(w)\tau^{\alpha-1}(w)\frac{\psi(w)-\varphi(w)}{1-\overline{\psi(w)}\varphi(w)}
-u(w)\varphi^\sharp(w)\frac{\psi(w)-\varphi(w)}{1-\overline{\psi(w)}\varphi(w)}\\
& &+v(w)\psi^\sharp(w)\frac{\psi(w)-\varphi(w)}{1-\overline{\varphi(w)}\psi(w)}
+u(w)\varphi^\sharp(w)\frac{\psi(w)-\varphi(w)}{1-\overline{\psi(w)}\varphi(w)}\Big|\\
&\geq&|u(w)\varphi^\sharp(w)+v(w)\psi^\sharp(w)|\rho(w)-|u(w)\varphi^\sharp(w)||1-\tau^{\alpha-1}(w)|\rho(w)\\
&\gtrsim&|u(w)\varphi^\sharp(w)+v(w)\psi^\sharp(w)|\rho(w)-C|u(w)\varphi^\sharp(w)|\rho^2(w).\\
\end {eqnarray*}
By Remark 2.1 and (a)of (i),we have that
\begin{eqnarray*}
|u(w)\varphi^\sharp(w)+v(w)\psi^\sharp(w)|\rho(w)\leq C+C|u(w)\varphi^\sharp(w)\rho^2(w)\leq C.
\end{eqnarray*}

{\bf Definition 3.1} Let$ \{z_n\}$ be a sequence in $\mathbb{D}$.

(i)Let $D_{u,\varphi}$ be the set of all sequences $ \{z_n\}$ such that $|u(z_n)\varphi^\sharp(z_n)|\rightarrow\infty$ when $|z_n|\rightarrow 1$.

(ii)Let $E_{u,\varphi}$ be the set of all sequences $ \{z_n\}$ such that
\begin{eqnarray*}
\frac{|U(z_n)|}{(1-|\varphi(z_n)|^2)^{\alpha-1}}\rightarrow\infty,
\end{eqnarray*}
when $|z_n|\rightarrow 1$.

{\bf Remark 3.1.} Using these notions, \cite{OSZ} is rewritten that a weighted composition operator $uC_\varphi$
is bounded from $\mathcal{B}^\alpha$ to $\mathcal{B}^\beta$ if and only if $D_{u,\varphi}=E_{u,\varphi}=\emptyset$.

We will characterize the boundedness of $uC_\varphi-vC_\psi$ from $\mathcal{B}^\alpha$ to $\mathcal{B}^\beta$ under the assumption
\begin{eqnarray*}
\sup\limits_{z\in\mathbb{D}}\frac{|U(z)|\rho(z)}{(1-|\varphi(z)|^2)^{\alpha-1}}<\infty\ \ \ \ \rm{and} \ \ \ \
\sup\limits_{z\in\mathbb{D}}\frac{|V(z)|\rho(z)}{(1-|\psi(z)|^2)^{\alpha-1}}<\infty.\ \ \ (\rm{A})
\end{eqnarray*}

{\bf Theorem 3.1}
Let $\varphi,\psi\in S(\mathbb{D})$ and $u,v\in H(\mathbb{D})$. Suppose that the condition (A) holds.
Then $uC_\varphi-vC_\psi$ is bounded  from $\mathcal{B}^\alpha(\alpha>1)$ to $\mathcal{B}^\beta$ if and only if
the following conditions hold:

(i)$D_{u,\varphi}=D_{v,\psi}$.

(ii)$\sup\limits_{z\in\mathbb{D}}|u(w)\varphi^\sharp(w)-v(w)\psi^\sharp(w)|<\infty$.

(iii)$(a)\sup\limits_{z\in\mathbb{D}}|u(w)\varphi^\sharp(w)|\rho(z)<\infty,\ \ \ \ (b)\sup\limits_{z\in\mathbb{D}}|v(w)\psi^\sharp(w)|\rho(z)<\infty$.

(iv)$E_{u,\varphi}=E_{v,\psi}$.

(v)
$
(a)
\sup\limits_{z\in E_{u,\varphi}}\lim\limits_{n\rightarrow\infty}\frac{|U(z_n)-V(z_n)|}{(1-|\varphi(z_n)|^2)^{\alpha-1}}<\infty,
(b)
\sup\limits_{z\in E_{v,\psi}}\lim\limits_{n\rightarrow\infty}\frac{|U(z_n)-V(z_n)|}{(1-|\psi(z_n)|^2)^{\alpha-1}}<\infty.
$

(vi)$(a)
\sup\limits_{z\in E_{u,\varphi}}\lim_{n\rightarrow\infty}|U(z_n)|d_\alpha(z_n)<\infty,
(b)
\sup\limits_{z\in  E_{v,\psi}}\lim_{n\rightarrow\infty}|V(z_n)|d_\alpha(z_n)<\infty.
$

{\bf Proof.} First we assume that $uC_\varphi-vC_\psi$ is bounded from  $\mathcal{B}^\alpha$ to $\mathcal{B}^\beta$.
In the following we prove (i)-(iii). Suppose that $D_{v,\psi}=\emptyset$. Then (iii)(b) holds. By (iii) of Proposition 3.2, (iii)(a) also holds.
Fix $w\in\mathbb{D}$, let $p(z)=f_\lambda(z)\sigma_\lambda(z)$, then $p\in\mathcal{B}^\alpha$. Take $\lambda=\varphi(w), z=w$.
By the estimate on $\|(uC_\varphi-vC_\psi)(p)\|_\beta$, we have that
\begin {eqnarray*}
C&\geq&\Big|u(w)\varphi^\sharp(w)+V(w)\sigma(w)\frac{(1-|\varphi(w)|^2)^{\alpha-1}}{(1-\overline{\varphi(w)}\psi(w))^{2(\alpha-1)}}
-v(w)\psi^\sharp(w)\tau^{\alpha}(w)\\
&-&2(\alpha-1)\overline{\varphi(w)}v(w)\psi^\sharp(w)\tau^{\alpha-1}(w)\sigma(w)
\frac{1-|\psi(w)|^2}{1-\overline{\varphi(w)}\psi(w)}\Big|.
\end {eqnarray*}
Hence
\begin {eqnarray} \label{3.1}
& &\Big|u(w)\varphi^\sharp(w)-v(w)\psi^\sharp(w)\tau^{\alpha}(w)\\  \notag\
&-&2(\alpha-1)\overline{\varphi(w)}v(w)\psi^\sharp(w)\tau^{\alpha-1}(w)\sigma(w)
\frac{1-|\psi(w)|^2}{1-\overline{\varphi(w)}\psi(w)}\Big|\\ \notag\
&\leq& C+|V(w)|\rho(w)\Big|\frac{(1-|\varphi(w)|^2)^{\alpha-1}}{(1-\overline{\varphi(w)}\psi(w))^{2(\alpha-1)}}\Big|.
\end {eqnarray}
With condition (A) and the assume of $|v(w)\psi^\sharp(w)|<\infty$, we have
\begin {eqnarray} \label{3.2}
|u(w)\varphi^\sharp(w)-v(w)\psi^\sharp(w)\tau^{\alpha-1}(w)|\leq C.
\end {eqnarray}
From Remark 2.1, we obtain
\begin{eqnarray}\label{3.3}
|u(w)\varphi^\sharp(w)-v(w)\psi^\sharp(w)|&\leq& C+|v(w)\psi^\sharp(w)(1-\tau^{\alpha-1}(w))|\leq C,
\end{eqnarray}
therefore
\begin {eqnarray*}
|u(w)\varphi^\sharp(w)|\leq|u(w)\varphi^\sharp(w)-v(w)\psi^\sharp(w)|+|v(w)\psi^\sharp(w)|
\leq C.
\end {eqnarray*}
Thus we get $G_{u,\varphi}=\emptyset$ which implies that $D_{v,\psi}\subset D_{u,\varphi}$. Similarly, we can obtain that
$D_{u,\varphi}\subset D_{v,\psi}$.
\

Next suppose that $D_{v,\psi}\neq\emptyset$.
Multiply inequality (3.1) by $\rho(w)$ and combine with (i) of Proposition 3.2, we have that
\begin {eqnarray*}
|u(w)\varphi^\sharp(w)-v(w)\psi^\sharp(w)\tau^\alpha(w)|\rho(w)\leq C.
\end {eqnarray*}
Therefore with Remak 2.1, we obtain
\begin {eqnarray*}
& &|u(w)\varphi^\sharp(w)-v(w)\psi^\sharp(w)|\rho(w)\\
&\leq&|u(w)\varphi^\sharp(w)-v(w)\psi^\sharp(w)\tau^{\alpha}(w)|\rho(w)+|v(w)\psi^\sharp(w)||1-\tau^{\alpha}(w)|\rho(w)\\
&\lesssim& C.
\end {eqnarray*}
Combining with (iii) of Proposition 3.2, we have that
$$
\sup\limits_{z\in\mathbb{D}}|u(w)\varphi^\sharp(w)|\rho(w)<\infty,\ \ \ \sup\limits_{z\in\mathbb{D}}|v(w)\psi^\sharp(w)|\rho(w)<\infty.
$$
Hence from inequality (3.1) again, we have that
\begin {eqnarray} \label{3.4}
|u(w)\varphi^\sharp(w)-v(w)\psi^\sharp(w)\tau^{\alpha}(w)|\leq C.
\end {eqnarray}
Applying Remark 2.1 again, we have
\begin {eqnarray*}
|u(w)\varphi^\sharp(w)-v(w)\psi^\sharp(w)|\leq C.
\end {eqnarray*}
Hence we have that $u(w)\varphi^\sharp(w)\rightarrow\infty$ which implies that $D_{v,\psi}\subset D_{u,\varphi}$. Similarly we can prove that $D_{u,\varphi}\subset D_{v,\psi}$. Hence we get the condition (i)-(iii).
\

In the following we prove (iv)-(vi).
For $f\in \mathcal{B}^\alpha$, by the triangle inequality, we have that
\begin {eqnarray*}
\|(uC_\varphi-vC_\psi)(f)\|_{\mathcal{B}^\beta}
&\geq&|U(z)f(\varphi(z))-V(z)f(\psi(z))|\\
&-&\|f\|_{\mathcal{B}^\alpha}|u(z)\varphi^\sharp(z)-v(z)\psi^\sharp(z)|.
\end {eqnarray*}
Hence
\begin {eqnarray}\label{3.5}
|U(z)f(\varphi(z))-V(z)f(\psi(z))|\leq \|(uC_\varphi-vC_\psi)(f)\|_{\mathcal{B}^\beta}+C\|f\|_{\mathcal{B}^\alpha}.
\end {eqnarray}
Consider the function $g_\lambda$ defined by $g_\lambda(z)=\frac{1}{(1-\overline{\lambda}z)^{\alpha-1}}$ and take $\lambda=\varphi(w), z=w$, then from the up inequality we have that
\begin {eqnarray}\label{3.6}
\Big|\frac{U(w)}{(1-|\varphi(w)|^2)^{\alpha-1}}
-\frac{V(w)}{(1-\overline{\varphi(w)}\psi(w))^{\alpha-1}}\Big|\leq C.
\end {eqnarray}
Consider the function $f_\lambda=\frac{(1-|\lambda|^2)^{\alpha-1}}{(1-\overline{\lambda}z)^{2(\alpha-1)}}$ and take $\lambda=\varphi(w),z=w$, similarly we have that
\begin {eqnarray}\label{3.7}
& &\Big|\frac{U(w)}{(1-|\varphi(w)|^2)^{\alpha-1}}
-V(w)\frac{(1-|\varphi(w)|^2)^{\alpha-1}}{(1-\overline{\varphi(w)}\psi(w))^{2(\alpha-1)}}\Big|\leq C.
\end {eqnarray}
By (\ref{3.6}) which multiplied with $\Big|\frac{1-|\varphi(w)|^2}{1-\overline{\varphi(w)}\psi(w)}\Big|^{\alpha-1}$ and (\ref{3.7}), we have that
\begin {eqnarray}\label{3.8}
|U(w)|\Big|\frac{1}{(1-\overline{\varphi(w)}\psi(w))^{\alpha-1}}-\frac{1}{(1-|\varphi(w)|^2)^{\alpha-1}}\Big|\leq C.
\end {eqnarray}
Similarly we have that
\begin {eqnarray}\label{3.9}
|V(w)|\Big|\frac{1}{(1-\overline{\psi(w)}\varphi(w))^{\alpha-1}}-\frac{1}{(1-|\psi(w)|^2)^{\alpha-1}}\Big|\leq C.
\end {eqnarray}
From inequalities (\ref{3.6}) and (\ref{3.8}), we have that
\begin {eqnarray}\label{3.10}
|U(w)-V(w)|\Big|\frac{1}{(1-\overline{\varphi(w)}\psi(w))^{\alpha-1}}\Big|\leq C.
\end {eqnarray}

First, suppose that $E_{u,\varphi}=\emptyset.$ Then by (\ref{3.10}) and Lemma 2.2, we obtain that
\begin {eqnarray*}
\Big|V(w)\frac{1}{(1-\overline{\varphi(w)}\psi(w))^{\alpha-1}}\Big|<\infty.
\end {eqnarray*}
Apply Lemma 2.2 again, we have that
\begin {eqnarray*}
\sup\limits_{z\in\mathbb{D}}\Big|V(z)\frac{1}{(1-|\psi(z)|^2)^{\alpha-1}}\Big|<\infty,
\end {eqnarray*}
which implies that $E_{v,\psi}=\emptyset.$
\

Next suppose that $E_{u,\varphi}\neq\emptyset$ and let $\{z_n\}\in E_{u,\varphi}$. By (\ref{3.6}), we get
\begin {eqnarray*}
|V(z_n)|\Big|\frac{1}{(1-\overline{\varphi(z_n)}\psi(z_n))^{\alpha-1}}\Big|\rightarrow\infty.
\end {eqnarray*}
Then from (\ref{3.9}), we get that
\begin {eqnarray*}
|V(z_n)|\Big|\frac{1}{(1-|\psi(z_n)|^2)^{\alpha-1}}\Big|\rightarrow\infty.
\end {eqnarray*}
This mean that $E_{u,\varphi}\subset E_{v,\psi}$. Similarly we obtain $E_{v,\psi}\subset E_{u,\varphi}$. Then we obtain the condition (iv).
\

Combining (\ref{3.6}) and (\ref{3.7}), for $\{z_n\}\in E_{u,\varphi}$,
\begin {eqnarray*}
\lim\limits_{n\rightarrow\infty}|V(z_n)|\Big|\frac{1}{(1-\overline{\varphi(z_n)}\psi(z_n))^{\alpha-1}}\Big|
\Big|1-\frac{(1-|\varphi(z_n)|^2)^{\alpha-1}}{(1-\overline{\varphi(z_n)}\psi(z_n))^{\alpha-1}}\Big|\leq C.
\end {eqnarray*}
Thus we get
\begin {eqnarray}\label{3.11}
\Big|\frac{(1-|\varphi(z_n)|^2)^{\alpha-1}}{(1-\overline{\varphi(z_n)}\psi(z_n))^{\alpha-1}}\Big|\rightarrow 1.
\end {eqnarray}
By (\ref{3.8}) and (\ref{3.11}), for all $\{z_n\}\in E_{u,\varphi}$, we have that
\begin {eqnarray*}
\lim\limits_{n\rightarrow\infty}|U(z_n)-V(z_n)|\frac{1}{(1-|\varphi(z_n)|^2)^{\alpha-1}}\leq C.
\end {eqnarray*}
Similarly we have that
\begin {eqnarray} \label{3.12}
\lim\limits_{n\rightarrow\infty}|U(z_n)-V(z_n)|\frac{1}{(1-|\psi(z_n)|^2)^{\alpha-1}}\leq C.
\end {eqnarray}
Thus we get the condition (v).

From (\ref{3.5}), we obtain that there exists a positive constant $C$ such that for any $\{z_n\}\in E_{u,\varphi}$,
and any $f \in \mathcal{B}^\alpha$ with $\|f\|_{\mathcal{B}^\alpha}\leq1$,
\begin {eqnarray*}
\lim\limits_{n\rightarrow\infty}|U(z_n)||f(\varphi(z_n))-f(\psi(z_n))|\leq C.
\end {eqnarray*}
Therefore we get
\begin {eqnarray*}
\sup\limits_{\{z_n\}\in E_{u,\varphi}}\lim\limits_{n\rightarrow\infty}|U(z_n)|d_\alpha(z_n)\leq C.
\end {eqnarray*}
Hence condition (vi) holds.
\

Conversely, we suppose that the condition (i)-(vi) hold. Let $\partial E_{u,\varphi}$
be the set of the cluster points of each $\{z_n\}\in E_{u,\varphi}$. Then $\partial E_{u,\varphi}\subset\partial\mathbb{D}$.
There exist a constant $M_1>0$ and a compact subset $K$ of $\overline{\mathbb{D}}$ such that $\partial E_{u,\varphi}\subset K$,
\begin {eqnarray*}
\sup\limits_{z\in K}|U(z)-V(z)|\frac{1}{(1-|\varphi(z)|^2)^{\alpha-1}}+\sup\limits_{z\in K}|V(z)|d_\alpha(z)\leq M_1,
\end {eqnarray*}
and
\begin {eqnarray*}
\sup\limits_{z\in K}|U(z)-V(z)|\frac{1}{(1-|\psi(z)|^2)^{\alpha-1}}+\sup\limits_{z\in K}|U(z)|d_\alpha(z)\leq M_1.
\end {eqnarray*}

On the other hand, there exits a constant $M_2>0$ such that
\begin {eqnarray*}
\sup\limits_{z\in\mathbb{D}\setminus K}|U(z)|\frac{1}{(1-|\varphi(z)|^2)^{\alpha-1}}
+\sup\limits_{z\in\mathbb{D}\setminus K}|V(z)|\frac{1}{(1-|\psi(z)|^2)^{\alpha-1}}\leq M_2.
\end {eqnarray*}

Then we obtain that for $M_0=\max\{M_1,M_2\}$,
\begin {eqnarray*}
\sup\limits_{z\in\mathbb{D}}L(z)<M_0.
\end {eqnarray*}

Hence, by proposition 3.1, $uC_\varphi-vC_\psi$ is bounded from $\mathcal{B}^\alpha$ to $\mathcal{B}^\beta$. This finishes the proof.

\section{The compactness from $\mathcal{B}^\alpha$ to $\mathcal{B}^\beta$}
In this section, we give the compactness of $uC_\varphi-vC_\psi$ from $\mathcal{B}^\alpha$ to $\mathcal{B}^\beta$. For this purpose, we need the following Lemma which can be proved in a standard way.(See, for example, Proposition 3.11 of \cite{CB})

{\bf Lemma 4.1} Let $\varphi,\psi$ be in $S(\mathbb{D})$ and $u,v$ be in $H(\mathbb{D})$. Suppose that $uC_\varphi-vC_\psi$ is bounded from $\mathcal{B}^\alpha$ to $\mathcal{B}^\beta$. Then the following are equivalent:

(i) $uC_\varphi-vC_\psi$ is compact from $\mathcal{B}^\alpha$ to $\mathcal{B}^\beta$.

(ii)$\|uC_\varphi-vC_\psi\|_{\mathcal{B}^\beta}\rightarrow 0$ for any bounded sequence $\{f_n\}$ in $\mathcal{B}^\alpha$
that converges to 0 uniformly on every compact subset of $\mathbb{D}$.

(iii)$\|uC_\varphi-vC_\psi\|_\beta\rightarrow 0$ for any sequence as in (ii).

We define some other sets of sequences in $\mathbb{D}$.

{\bf Definition 4.2.} Let $\varphi$ be in $S(\mathbb{D})$ and $u$ be in $H(\mathbb{D})$. Let $\{z_n\}$ be a sequence in $\mathbb{D}$.

(i) We denote by $G_{u,\varphi}$ the set of all sequence $\{z_n\}$ such that $|\varphi(z_n)|\rightarrow 1$ and
$
|u(z_n)\varphi^\sharp(z_n)|\nrightarrow 0.
$

(ii) We denote by $F_{u,\varphi}$ the set of all $\{z_n\}$ such that $|\varphi(z_n)|\rightarrow 1$ and\\
$
|u(z_n)|\frac{1}{(1-|\varphi(z_n)|^2)^{\alpha-1}}\nrightarrow 0.
$

{\bf Remark 4.1.} The weighted operator $uC_\varphi$ is compact from $\mathcal{B}^\alpha$ to $\mathcal{B}^\beta$ if and only if $G_{u,\varphi}=F_{u,\varphi}=\emptyset$.
\

To characterize the compactness of $uC_\varphi-vC_\psi$  from $\mathcal{B}^\alpha$ to $\mathcal{B}^\beta$, we need the assumpation
\begin {eqnarray*}
\frac{|U(z_n)|\rho_n}{(1-|\varphi_n|^2)^{\alpha-1}}\rightarrow 0\ \ \ \  \rm{and} \ \ \ \
\frac{|V(z_n)|\rho_n}{(1-|\psi_n|^2)^{\alpha-1}}\rightarrow 0\ \ \ \ (\rm{B})
\end {eqnarray*}
 when $\rho_n\rightarrow 0$.


{\bf Theorem 4.3.} Let $\varphi,\psi$ be in $S(\mathbb{D})$ and $u,v$ be in $H(\mathbb{D})$ satisfying the condition (B). Suppose that $uC_\varphi-vC_\psi$ is bounded from $\mathcal{B}^\alpha(\alpha>1)$ to $\mathcal{B}^\beta$. Then $uC_\varphi-vC_\psi$ is compact from $\mathcal{B}^\alpha$ to $\mathcal{B}^\beta$ if and only if the following conditions hold:\\
(i)$G_{u,\varphi}=G_{v,\psi}$.\\
(ii) For any $\{z_n\}\in G_{u,\varphi}$,
\begin {eqnarray*}
\lim\limits_{n\rightarrow\infty}|u(z_n)\varphi^\sharp(z_n)-v(z_n)\psi^\sharp(z_n)|=0.
\end {eqnarray*}
(iii) For any $\{z_n\}\in G_{u,\varphi}$,
\begin {eqnarray*}
(a)\lim\limits_{n\rightarrow\infty}|u(z_n)\varphi^\sharp(z_n)|\rho_n=0,\ \
(b)\lim\limits_{n\rightarrow\infty}|v(z_n)\psi^\sharp(z_n)|\rho_n=0.
\end {eqnarray*}
(iv))$F_{u,\varphi}=F_{v,\psi}$.\\
(v) For any $\{z_n\}\in F_{u,\varphi}$,
\begin {eqnarray*}
(a)\lim\limits_{n\rightarrow\infty}\frac{|U(z_n)-V(z_n)|}{(1-|\varphi_n|^2)^{\alpha-1}}=0,\ \
(b)\lim\limits_{n\rightarrow\infty}\frac{|U(z_n)-V(z_n)|}{(1-|\psi_n|^2)^{\alpha-1}}=0.
\end {eqnarray*}
(vi) For any $\{z_n\}\in F_{u,\varphi}$,
\begin {eqnarray*}
(a)\lim\limits_{n\rightarrow\infty}|U(z_n)|d_\alpha(z_n)=0,\ \
(b)\lim\limits_{n\rightarrow\infty}|V(z_n)|d_\alpha(z_n)=0.
\end {eqnarray*}

{\bf Proof.} Suppose that $uC_\varphi-vC_\psi$ is compact from $\mathcal{B}^\alpha$ to $\mathcal{B}^\beta$. Consider the function defined by $f_{m,k}(z)=f_\lambda(z)(\sigma_\lambda^m(z)-\lambda^{m-k}\sigma_\lambda^k(z))$, for $\lambda\in \mathbb{D}$ and non-negative integers $k,m$ with $m>k\geq 1$. It is easy to see $f_{m,k}\in\mathcal{B}^\alpha$ converges to 0 uniformly on every compact subset of $\mathbb{D}$. Take $\lambda=\varphi(z_n), z=z_n, \{z_n\}\in\mathbb{D}$, then we have that
\begin {eqnarray*}
& &\|(uC_\varphi-vC_\psi)(f_{m,k})(z_n)\|_\beta\\
&=&\Big|ku(z_n)\varphi^\sharp(z_n)\varphi^{m-k}_n\sigma_{\varphi_n}^{k-1}(\varphi_n)
-V(z_n)\frac{(1-|\varphi_n|^2)^{\alpha-1}}{(1-\overline{\varphi_n}\psi_n)^{2(\alpha-1)}}\sigma_n^k[\sigma_n^{m-k}-\varphi_n^{m-k}]\\
& &-2(\alpha-1)\overline{\varphi_n}v(z_n)\psi^\sharp(z_n)
\frac{(1-|\varphi_n|^2)^{\alpha-1}(1-|\psi_n|^2)^\alpha}{(1-\overline{\varphi_n}\psi_n)^{2\alpha-1}}
\sigma_n^k[\sigma_n^{m-k}-\varphi_n^{m-k}]\\
& &+v(z_n)\psi^\sharp(z_n)\sigma_n^{k-1}\tau_n^\alpha[m\sigma_n^{m-k}-k\varphi_n^{m-k}]\Big|.
\end {eqnarray*}
Choose $m=4$ and $k=3$, we have that
\begin {eqnarray}\label{4.1}
& &|v(z_n)\psi^\sharp(z_n)\sigma_n^2[4\sigma_n-3\varphi_n]\tau_n^\alpha-V(z_n)
\frac{(1-|\varphi_n|^2)^{\alpha-1}}{(1-\overline{\varphi_n}\psi_n)^{2(\alpha-1)}}\sigma^3_n[\sigma_n-\varphi_n]\\ \notag\
& &-v(z_n)\psi^\sharp(z_n)\sigma_n^32(\alpha-1)\overline{\varphi_n}
\frac{(1-|\varphi_n|^2)^{\alpha-1}(1-|\psi_n|^2)^\alpha}{(1-\overline{\varphi_n}\psi_n)^{2\alpha-1}}|\rightarrow 0.
\end {eqnarray}
Choose $m=3$ and $k=2$, we have that
\begin {eqnarray}\label{4.2}
& &|v(z_n)\psi^\sharp(z_n)\sigma_n[3\sigma_n-2\varphi_n]\tau_n^\alpha-V(z_n)
\frac{(1-|\varphi_n|^2)^{\alpha-1}}{(1-\overline{\varphi_n}\psi_n)^{2(\alpha-1)}}\sigma^2_n[\sigma_n-\varphi_n]\\ \notag\
& &-v(z_n)\psi^\sharp(z_n)\sigma_n^22(\alpha-1)\overline{\varphi_n}
\frac{(1-|\varphi_n|^2)^{\alpha-1}(1-|\psi_n|^2)^\alpha}{(1-\overline{\varphi_n}\psi_n)^{2\alpha-1}}|\rightarrow 0.
\end {eqnarray}
Formula (\ref{4.1}) minus formula (\ref{4.2}) which multiplied by $\sigma_n$, then we have that
\begin {eqnarray}\label{4.3}
|v(z_n)\psi^\sharp(z_n)\sigma_n^2(\sigma_n-\varphi_n)\tau_n^\alpha|\rightarrow0.
\end {eqnarray}
Formula (\ref{4.2}) multiplied by $\sigma_n$ minus (\ref{4.3}) multiplied by 3, then we have that
\begin {eqnarray}\label{4.4}
& &|v(z_n)\psi^\sharp(z_n)\sigma_n^2\varphi_n\tau_n^\alpha-V(z_n)\frac{(1-|\varphi_n|^2)^{\alpha-1}}{(1-\overline{\varphi_n}\psi_n)^{2(\alpha-1)}}
\sigma^3_n(\sigma_n-\varphi_n)\\ \notag\
& &-v(z_n)\psi^\sharp(z_n)\sigma_n^3(\sigma_n-\varphi_n)2(\alpha-1)\overline{\varphi_n}
\frac{(1-|\varphi_n|^2)^{\alpha-1}(1-|\psi_n|^2)^\alpha}{(1-\overline{\varphi_n}\psi_n)^{2\alpha-1}}|\rightarrow 0.
\end {eqnarray}
Choose $m=2$ and $k=1$, we obtain that
\begin {eqnarray}\label{4.5}
& &|v(z_n)\psi^\sharp(z_n)(2\sigma_n-\varphi_n)\tau_n^\alpha-V(z_n)
\frac{(1-|\varphi_n|^2)^{\alpha-1}}{(1-\overline{\varphi_n}\psi_n)^{2(\alpha-1)}}
\sigma_n(\sigma_n-\varphi_n)\\ \notag\
& &-v(z_n)\psi^\sharp(z_n)\sigma_n(\sigma_n-\varphi_n)2(\alpha-1)\overline{\varphi_n}
\frac{(1-|\varphi_n|^2)^{\alpha-1}(1-|\psi_n|^2)^\alpha}{(1-\overline{\varphi_n}\psi_n)^{2\alpha-1}}\\  \notag\
& &+u(z_n)\varphi^\sharp(z_n)\varphi_n|\rightarrow0.
\end {eqnarray}
Formula (\ref{4.1}) minus formula (\ref{4.5}) which multiplied by $\sigma^2_n$, then we have that
\begin {eqnarray}\label{4.6}
|v(z_n)\psi^\sharp(z_n)\sigma_n^2(2\sigma_n-2\varphi_n)\tau_n^\alpha-u(z_n)\varphi^\sharp(z_n)\varphi_n\sigma_n^2|\rightarrow0.
\end {eqnarray}
From formulas (\ref{4.3}) and (\ref{4.6}), we obtain that
\begin {eqnarray}\label{4.7}
|u(z_n)\varphi^\sharp(z_n)\sigma_n^2|=|u(z_n)\varphi^\sharp(z_n)|\rho_n^2\rightarrow0.
\end {eqnarray}

Suppose that $G_{u,\varphi}\neq\emptyset$ and let $\{z_n\}\in G_{u,\varphi}$. Then (\ref{4.7}) implies that $\rho_n\rightarrow0$ or $|\tau_n|=1-|\rho^2_n|\rightarrow1$. In fact, the function $p(z)=f_\lambda(z)\sigma_\lambda(z)$ defined in Theorem 3.1 converges to 0 uniformly on every compact subset of $\mathbb{D}$ as $|\lambda|\rightarrow1$. Take $\lambda=\varphi(z_n), z=z_n, $ hence $0\leftarrow\|(uC_\varphi-vC_\psi)(p)\|_{\mathcal{B}^\beta}$ which implies that
\begin {eqnarray}\label{4.8}
0&\leftarrow&\|(uC_\varphi-vC_\psi)(p)\|_{\mathcal{B}^\beta}+|V(z_n)|\rho_n
\Big|\frac{(1-|\varphi_n|^2)^{\alpha-1}}{(1-\overline{\varphi_n}\psi_n)^{2(\alpha-1)}}\Big|\\ \notag\
&\geq&|u(z_n)\varphi^\sharp(z_n)-v(z_n)\psi^\sharp(z_n)\tau^{\alpha}_n|\\ \notag\
&-&\Big|2(\alpha-1)\overline{\varphi_n}v(z_n)\psi^\sharp(z_n)\tau^{\alpha-1}_n\sigma_n\frac{1-|\psi_n|^2}{1-\overline{\varphi_n}\psi_n}\Big|.\\ \notag\
\end {eqnarray}
Multiplying (\ref{4.8}) by $\rho^2_n$ and combining with (i)b of Proposition 3.2, we have that
\begin {eqnarray}\label{4.9}
|u(z_n)\varphi^\sharp(z_n)-v(z_n)\psi^\sharp(z_n)\tau^{\alpha}_n|\rho^2_n\rightarrow0.
\end {eqnarray}
Hence
\begin {eqnarray}\label{4.10}
|v(w)\psi^\sharp(w)\tau^{\alpha}(w)|\rho^2_n&\leq&
|u(w)\varphi^\sharp(w)-v(w)\psi^\sharp(w)\tau^{\alpha}(w)|\rho^2_n\\
&+&|u(w)\varphi^\sharp(w)|\rho^2_n \notag\
\rightarrow0,
\end {eqnarray}
which implies that
\begin {eqnarray} \label{4.11}
|v(w)\psi^\sharp(w)|\rho^2_n\rightarrow0.
\end {eqnarray}
Multiplying (\ref{4.8}) by $\rho_n$ and combining with (\ref{4.11}), we obtain
\begin {eqnarray}\label{4.12}
|u(z_n)\varphi^\sharp(z_n)-v(z_n)\psi^\sharp(z_n)\tau^{\alpha}_n|\rho_n\rightarrow0.
\end {eqnarray}
 From condition (B) and (\ref{4.4}), we obtain that
\begin {eqnarray}\label{4.13}
|v(z_n)\psi^\sharp(z_n)|\rho_n\rightarrow0.
\end {eqnarray}
Hence, applying (\ref{4.8}) again, we have that
\begin {eqnarray}\label{4.14}
|u(z_n)\varphi^\sharp(z_n)-v(z_n)\psi^\sharp(z_n)\tau^{\alpha}_n|\rightarrow0.
\end {eqnarray}
Therefore
\begin {eqnarray*}
0&\leftarrow&|u(z_n)\varphi^\sharp(z_n)-v(z_n)\psi^\sharp(z_n)\tau_n^\alpha|\\
&\geq&|u(z_n)\varphi^\sharp(z_n)-v(z_n)\psi^\sharp(z_n)|-|v(z_n)\psi^\sharp(z_n)||1-\tau_n^\alpha|\\
&\gtrsim&|u(z_n)\varphi^\sharp(z_n)-v(z_n)\psi^\sharp(z_n)|-C|v(z_n)\psi^\sharp(z_n)|\rho_n,
\end {eqnarray*}
which implies that
\begin {eqnarray}\label{4.15}
|u(z_n)\varphi^\sharp(z_n)-v(z_n)\psi^\sharp(z_n)|\rightarrow0.
\end {eqnarray}
Hence we get (ii) and (iii). Equation (\ref{4.15}) means that $\emptyset\neq G_{u,\varphi}\subset G_{v,\psi}$. Since the converse inclusion is also true,
$G_{u,\varphi}=G_{v,\psi}$ if $G_{u,\varphi}\neq\emptyset$.
If $G_{u,\varphi}=\emptyset$, then the inclusion implies that $G_{v,\psi}=\emptyset$ which proof is similar to Theorem 3.1.

In the following, we prove (iv)-(vi). Similar to Theorem 3.1, for any bounded sequence $\{g_n\}\in\mathcal{B}^\alpha$ that converges to 0 uniformly on every compact subset of $\mathbb{D}$, we have that
\begin {eqnarray}\label{4.16}
|U(z_n)g_n(\varphi(z))-V(z_n)g_n(\psi(z))|\rightarrow0.
\end {eqnarray}
Let $\{z_n\}\in\mathbb{D}$ such that $|\varphi(z_n)|\rightarrow1$ and let
\begin {eqnarray*}
f_k(z)=(\frac{1}{(1-|\varphi(z_n)|^2)^{\alpha-1}})^{-k+1}(\frac{1}{(1-\overline{\varphi(z_n)}z)^{\alpha-1}})^k
\end {eqnarray*}
for any positive integer $k$. Fix $k$, $\{f_k\}$ is bounded in $\mathcal{B}^\alpha$ and converges to 0 uniformly on every compact subsets of $\mathbb{D}$ as $n\rightarrow\infty$. For $z=z_n$ and $k=1,2$, we have that
\begin {eqnarray}\label{4.17}
\Big|\frac{U(z_n)}{(1-|\varphi(z_n)|^2)^{\alpha-1}}
-\frac{V(z_n)}{(1-\overline{\varphi(z_n)}\psi(z_n))^{\alpha-1}}\Big|\rightarrow0,
\end {eqnarray}
and
\begin {eqnarray}\label{4.18}
\Big|\frac{U(z_n)}{(1-|\varphi(z_n)|^2)^{\alpha-1}}
-\frac{V(z_n)}{(1-\overline{\varphi(z_n)}\psi(z_n))^{\alpha-1}}\frac{(1-|\varphi(z_n)|^2)^{\alpha-1}}{(1-\overline{\varphi(z_n)}\psi(z_n))^{\alpha-1}}\Big|
\rightarrow0.
\end {eqnarray}
The above formulas imply that
\begin {eqnarray}\label{4.19}
\Big|\frac{(1-|\varphi(z_n)|^2)^{\alpha-1}}{(1-\overline{\varphi(z_n)}\psi(z_n))^{\alpha-1}}\Big|\rightarrow1.
\end {eqnarray}
From (\ref{4.17}) and (\ref{4.18}), we have obtain that
\begin {eqnarray}\label{4.20}
|U(z_n)|\Big|\frac{1}{(1-\overline{\varphi(z_n)}\psi(z_n))^{\alpha-1}}-\frac{1}{(1-|\varphi(z_n)|^2)^{\alpha-1}}\Big|\rightarrow0.
\end {eqnarray}
Similarly we have that
\begin {eqnarray}\label{4.21}
|V(z_n)|\Big|\frac{1}{(1-\overline{\psi(z_n)}\varphi(z_n))^{\alpha-1}}-\frac{1}{(1-|\psi(z_n)|^2)^{\alpha-1}}\Big|\rightarrow0.
\end {eqnarray}
From (\ref{4.18}) we have that
\begin {eqnarray*}
\frac{|U(z_n)-V(z_n)|}{(1-|\varphi(z_n)|^2)^{\alpha-1}}\rightarrow0.
\end {eqnarray*}
Suppose that $\{z_n\}\in F_{u,\varphi}$. Then by (\ref{4.17}), we obtain that
\begin {eqnarray*}
\Big|\frac{V(z_n)}{(1-|\psi(z_n)|^2)^{\alpha-1}}\Big|
=\Big|\frac{V(z_n)}{(1-\overline{\varphi(z_n)}\psi(z_n))^{\alpha-1}}
\frac{(1-|\varphi(z_n)|^2)^{\alpha-1}}{(1-\overline{\varphi(z_n)}\psi(z_n))^{\alpha-1}}\Big|\rightarrow0.
\end {eqnarray*}
Hence, we get $\{z_n\}\in F_{v,\psi}$ which implies $F_{u,\varphi}\subset F_{v,\psi}$

Similar to the proof of Theorem 3.1, we get (iv) and (v)(b). Hence, by v(a) and the triangle inequality, we get (vi).

Conversely assume that the condition (i)-(vi) hold. Let $\{f_n\}$ be a sequence in $\mathcal{B}^\alpha$ such that $\|f_n\|_{\mathcal{B}^\alpha}\leq1$ and converges to 0 uniformly on every compact subsets of $\mathbb{D}$. Suppose $\|uC_\varphi-vC_\psi\|_{\mathcal{B}^\alpha}\nrightarrow0$. Then there exists a constant $\varepsilon_0>0$ such that  $\|uC_\varphi-vC_\psi(f_n)\|_{\mathcal{B}^\alpha}\geq\varepsilon_0$ for all $n$. Choose $z_n\in\mathbb{D}$ such that
\begin {eqnarray}\label{4.22}
(1-|z_n|^2)^\beta|(uf_n(\varphi_n)-vf_n(\psi_n))'(z_n)|\geq\varepsilon_0
\end {eqnarray}
for each $n$.
By condition (i)-(iii),
\begin {eqnarray*}
& &|u(z_n)\varphi^\sharp(z_n)(1-|\varphi(z_n)|^2)^\alpha f'_n(\varphi(z_n))-v(z_n)\psi^\sharp(z_n)(1-|\psi(z_n)|^2)^\alpha f'_n(\psi(z_n))|\\
&\lesssim&|u(z_n)\varphi^\sharp(z_n)-v(z_n)\psi^\sharp(z_n)|\|f_n\|_{\mathcal{B}^\alpha}+C|v(z_n)\psi^\sharp(z_n)|\rho(z_n)\\
&\rightarrow&0.
\end {eqnarray*}
By condition (iv)-(vi)
\begin {eqnarray*}
& &|U(z_n)f_n(\varphi(z_n))-V(z_n)f_n(\psi(z_n))|\\
&\leq&|U(z_n)|d_\alpha(z_n)+|U(z_n)-V(z_n)||f_n(\psi(z_n))|\\
&\rightarrow&0.
\end {eqnarray*}
Hence, we have that
\begin {eqnarray*}
& &(1-|z_n|^2)^\beta|(uf_n(\varphi)-vf_n(\psi))'(z_n)|\leq|U(z_n)f_n(\varphi(z_n))-V(z_n)f_n(\psi(z_n))|\\
&+&|u(z_n)\varphi^\sharp(z_n)(1-|\varphi(z_n)|^2)^\alpha f'_n(\varphi(z_n))-v(z_n)\psi^\sharp(z_n)(1-|\psi(z_n)|^2)^\alpha f'_n(\psi(z_n))|\\
&\rightarrow&0,
\end {eqnarray*}
which contradicts to (\ref{4.22}). Our proof is finished.


\end{document}